\documentclass[12pt,a4paper,openright,reqno]{extarticle}
\usepackage[english]{babel}
\usepackage{newlfont}
\usepackage{amsbsy}
\usepackage[dvips]{graphicx}
\usepackage{color}
\usepackage{bbm}
\usepackage{graphicx, subfigure}
\usepackage{color}
\usepackage[center]{caption2}
\usepackage{amssymb}
\usepackage{amsmath}

\usepackage{latexsym}
\usepackage{amsthm}
\usepackage{amsthm}

\usepackage[dvips]{geometry}
\theoremstyle{plain}
\newtheorem{thm}{Theorem}
\newtheorem{cor}[thm]{Corollary}

\newtheorem{lem}[thm]{Lemma}

\newtheorem{nota}[thm]{Notation}

\theoremstyle{remark}

\newcommand{\R}{\mathbb{R}}
\newcommand{\N}{\mathbb{N}}

\newcommand{\bp}{\begin{pmatrix}}
\newcommand{\ep}{\end{pmatrix}}

\usepackage{mathtools}
\def\alphaltiset#1#2{\ensuremath{\left(\kern-.2em\left(\genfrac{}{}{0pt}{}{#1}{#2}\right)\kern-.2em\right)}}
\usepackage[center]{caption2}

\begin{document}

\title{Affine subspaces of antisymmetric matrices with constant rank }
\author{Elena Rubei}
\date{}
\maketitle

{\footnotesize\em Dipartimento di Matematica e Informatica ``U. Dini'', 
viale Morgagni 67/A,
50134  Firenze, Italia }

{\footnotesize\em
E-mail address: elena.rubei@unifi.it}

\def\thefootnote{}
\footnotetext{ \hspace*{-0.36cm}


{\bf \copyright}   This is an Accepted Manuscript version of the following article, published in Linear and Multilinear Algebra:  E. Rubei "Affine subspaces of antisymmetric matrices with constant rank" Linear and Multilinear Algebra 72 (11), 1741-1750 (2024), DOI: 10.1080/03081087.2023.2198759. It is deposited under the terms of the Creative Commons Attribution-NonCommercial-NoDerivatives License (http://creativecommons.org/licenses/by-nc-nd/4.0/), which permits non-commercial re-use, distribution, and reproduction in any medium, provided the original work is properly cited, and is not altered, transformed, or built upon in any way}

\begin{abstract} 
For every $n \in \mathbb{N}$ and every field $K$,  let $A(n,K)$ be the vector space of the antisymmetric $(n \times n)$-matrices over  $K$.
We say that an affine subspace $S$ of  $A(n,K)$ has constant rank $r$ if every matrix of $S$ has rank $r$.
Define $${\cal A}_{antisym}^K(n;r)=  \{ S \;| \; S \; \mbox{\rm  affine subspace of $A(n,K)$  of constant rank }  r\}$$
$$a_{antisym}^K(n;r) = \max \{\dim S \mid S \in {\cal A}_{antisym}^K(n;r)  \}.$$

In this paper we prove the following formulas:

 for $n \geq 2r +2 $
$$a_{antisym}^{\mathbb{R}}( n; 2r) = (n-r-1) r   ;$$ 

for $n=2r$ 
$$a_{antisym}^{\mathbb{R}}( n; 2r) =r(r-1)   ;$$
 
for $n=2r+1$ 
$$a_{antisym}^{\mathbb{R}}( n; 2r) = r(r+1)   .$$ 

\end{abstract}

\section{Introduction}

For every $m,n \in \mathbb{N}$ and every field $K$, let $M(m \times n, K)$ be the vector space of the $(m \times n)$-matrices over $K$, let $S(n,K)$ be the vector space of the symmetric $(n \times n)$-matrices over  $K$ and let $A(n,K)$ be the vector space of the antisymmetric $(n \times n)$-matrices over  $K$. Moreover, denote the $\mathbb{R}$-vector space of the hermitian $(n\times n)$-matrices by $H(n)$.

We say that an affine subspace $S$ of $M(m \times n, K)$ has constant rank $r$ if every matrix of $S$ has rank $r$ and we say that a linear subspace $S$ 
 of $M(m \times n, K)$ has constant rank $r$ if every nonzero matrix of $S$ has rank $r$.

Define
$${\cal A}^K(m \times n; r)=  \{ S \;| \; S \; \mbox{\rm  affine subsapce of $M(m \times n, K)$  of constant rank }  r\}$$
$${\cal A}_{sym}^K(n;r)=  \{ S \;| \; S \; \mbox{\rm  affine subsapce of $S(n,K)$  of constant rank }  r\}$$
$${\cal A}_{antisym}^K(n;r)=  \{ S \;| \; S \; \mbox{\rm  affine subsapce of $A(n,K)$  of constant rank }  r\}$$
$${\cal A}_{herm}(n;r)=  \{ S \;| \; S \; \mbox{\rm  affine subsapce of $H(n)$  of constant rank }  r\}$$
$${\cal A}_{sym}^{\mathbb{R}}(n;p, \nu)=  \{ S \;| \; S \; \mbox{\rm  affine subsapce of $S(n,\mathbb{R})$ s.t.
each $ A \in S$  has signature }  (p, \nu)\}$$
$${\cal A}_{herm}(n;p, \nu)=  \{ S \;| \; S \; \mbox{\rm  affine subsapce of $H(n)$ s.t.
each $ A \in S$  has signature }  (p, \nu)\}$$
$${\cal L}^K (m \times n;r) =  \{ S \;| \; S \; \mbox{\rm  linear subsapce of $M(m \times n, K)$   of constant rank }  r\}$$
$${\cal L}_{sym}^K (n; r) =  \{ S \;| \; S \; \mbox{\rm  linear subsapce of $S(n,K)$  of constant rank }  r\}$$

Let 
$$a^K(m \times n;r) = \max \{\dim S \mid S \in {\cal A}^K(m \times n; r ) \}$$
$$a_{sym}^K(n;r) = \max \{\dim S \mid S \in {\cal A}_{sym}^K(n;r)  \}$$
$$a_{antisym}^K(n;r) = \max \{\dim S \mid S \in {\cal A}_{antisym}^K(n,r)  \}.$$
$$a_{herm}(n;r) = \max \{\dim S \mid S \in {\cal A}_{herm}(n;r)  \}$$
$$a_{sym}^{\mathbb{R}}(n;p, \nu) = \max \{\dim S \mid S \in {\cal A}_{sym}^{\mathbb{R}}(n; p,\nu)  \}$$
$$a_{herm}(n;p, \nu) = \max \{\dim S \mid S \in {\cal A}_{herm}(n; p,\nu)  \}$$
$$l^K(m \times n;r) = \max \{\dim S \mid S \in {\cal L}^K (m \times n;r)\}$$
$$l_{sym}^K(n;r) = \max \{\dim S \mid S \in {\cal L}_{sym}^K (n,r)\}.$$

There is a wide literature on linear subspaces of constant rank. In particular we quote the following theorems:

\begin{thm}  {\bf (Westwick, \cite{W1})} For $2 \leq r \leq m \leq n$, we have:
$$ n-r+1 \leq  l^{\mathbb{C}}(m \times n;r) \leq  m+ n -2 r+1$$
\end{thm}

\begin{thm}  {\bf (Ilic-Landsberg, \cite{I-L})} If $r$ is even and greater than or equal to $2$, then
$$l_{sym}^{\mathbb{C}}(n;r) =n-r +1$$
\end{thm}

In case $r$ odd, the following result holds, see \cite{I-L}, \cite{G}, \cite{H-P}:

\begin{thm}  If $r$ is odd, then
$$l_{sym}^{\mathbb{C}}(n;r) =1$$
\end{thm}

We mention also that,  in \cite{Fl}, Flanders
proved that, if $ r \leq m \leq n$, a linear subspace of $M(m \times n ,\mathbb{C})$ such that every of its elements has rank less than or equal to $r$ has dimension less than or equal to $r n$.

In \cite{Ru} we proved the following theorems:

\begin{thm}\label{mio}  Let $n,r \in \mathbb{N}$ with $r \leq n$.
Then
$$a_{sym}^{\mathbb{R}}(n;r) \leq
\left\lfloor \frac{r}{2} \right\rfloor
\left(n-  \left\lfloor \frac{r}{2} \right\rfloor\right)
.$$ 
\end{thm}

\begin{thm}\label{mio2}  Let $m,n,r \in \mathbb{N}$ with $r \leq m \leq n$.
Then
$$a^{\mathbb{R}}(m \times n;r) = rn- \frac{r(r+1)}{2} .$$ 
\end{thm}

We proved also a statement 
on  the maximal dimension of affine subspaces with constant signature in $S(n,\R)$ and one 
on the maximal dimension of affine  subspaces of constant rank in $H(n)$.

In this paper we investigate on the maximal dimension of affine subspaces in the space of the antisymmetric matrices of constant rank.  
Precisely we prove the following theorem:

\begin{thm} \label{thmantisym}
For $n \geq 2r+2  $: 
$$a_{antisym}^{\mathbb{R}}( n; 2r) = (n-r-1) r   .$$ 
For $n=2r$ 
$$a_{antisym}^{\mathbb{R}}( n; 2r) =r(r-1)   .$$ 
For $n=2r+1$ 
$$a_{antisym}^{\mathbb{R}}( n; 2r) = r(r+1)   .$$ 
\end{thm}

\section{Proof of the theorem}

\begin{nota} Let $m,n \in \mathbb{N} -\{0\} $. 

We denote the $n \times n$ identity matrix over $\R $ by $I_n$.

We denote $E_{i,j}^n$ the $n \times n$  matrix over $\R$ such that 
$$ (E_{i,j}^{n})_{x,y} = \left\{ \begin{array}{ll}
1 & \mbox{\rm if} \; (x,y)=(i,j) \\
0 & \mbox{\rm otherwise}
\end{array} \right.$$
We omit the superscript when it is clear from the context.

 We denote by $0_{m \times n }$ the $m \times n$ null real matrix. 
We omit the subscript when it is clear from the context.
 
 We denote $diag (d_1, \ldots, d_n)$ the diagonal matrix whose diagonal entries are $d_1, \ldots, d_n$.
 
For any $A \in M(m \times n, \R) $ we denote the submatrix of $A$ given by the rows $i_1, , \ldots , i_k$ and the columns $j_1, \ldots, j_s$  by $A^{(j_1, \ldots, j_s)}_{(i_1, \ldots, i_k)}$.

We denote by $J$ the $(2 \times 2)$-matrix 
$$ \bp 0 & 1 \\ -1 & 0 \ep $$
and we denote by $\overline{J}_{2n}$ the $(2n \times 2n)$ block diagonal matrix whose diagonal blocks are equal to $J$.
We omit the subscript when it is clear from the context.

We say that a  $(2 \times 2)$-matrix is {\em pinco} if and only if it is equal to $$ \bp
a & b \\ b & -a
 \ep $$
 for some $a,b \in \R$.
For any pinco  matrix $ P=\bp
a & b \\ b & -a
 \ep $, we denote by $\tilde{P}$ the matrix $JP$, that is the matrix 
 $$ \bp
b & -a \\ -a & -b
 \ep.$$
 
We say that a  $(2 \times 2)$-matrix is {\em antipinco} if and only if it is equal to $$ \bp
a & b \\ -b & a
 \ep $$
 for some $a,b \in \R$.

Let  $A_{i,j}$ for $i,j \in \{1, \ldots, n\}$ with $i <j$ be 
$2 \times 2$ matrices. We denote by $R(A_{i,j})$ the antisymmetric $(2n \times 2n )$-matrix 
$$ \bp
 0 & A_{1,2} & A_{1,3} & . & . & . \\ 
  -{}^t \!A_{1,2} & 0  & A_{2,3} & . & . & . \\ 
  -{}^t \! A_{1,3} &  -{}^t \!A_{2,3} & 0  & . & . & . \\ 
  . & . & . &  .& . &.  \\
   . & . & . & . & . &. \\
    . & . & . & . & . & .
 \ep.$$
 
 We denote by $X(A_{i,j})$ the symmetric $(2n \times 2n )$-matrix 
$$ \bp 0 & A_{1,2} & A_{1,3} & . & . & . \\  {}^t \!A_{1,2} & 0  & A_{2,3} & . & . & . \\  {}^t \! A_{1,3} &  {}^t \!A_{2,3} & 0  & . & . & . \\ . & . & . &  .& . &.  \\ . & . & . & . & . &. \\  . & . & . & . & . & . \ep.$$

 For any $l_1, \ldots,l_{r} \in \R$, let $T(l_1, \ldots,l_r)$ be the $2r \times 2r $ diagonal block matrix 
 whose diagnal blocks are 
 $l_1 J$,...., $l_rJ$.
\end{nota}



\begin{lem} \label{lemma2}
Let $m,r \in \mathbb{N}-\{0\}$.
Let $V$ be vector subspace of $\R^{3m}$.

Let $\pi_1: \R^{3m} \rightarrow \R^{2m}$ be the projection onto the first 
$2m$ coordinates.

Let $\pi_2: \R^{3m} \rightarrow \R^{2m} $ be the projection onto the last 
$2m$ coordinates.

Let $\pi_3 :  \R^{3m} \rightarrow \R^{2m}$ be the projection onto the  coordinates $1, \ldots, m, 2m+1, \ldots 3m$.

If $\dim (\pi_i (V) ) \leq 2r $ for $i=1,2,3$, then $\dim(V) \leq 3r$.
\end{lem}

\begin{proof}
Suppose to the contrary that $\dim(V) \geq 3r+1$. Then we would have:
$$\dim (Ker(\pi_i|_V))= \dim(V)- \dim (Im(\pi_i|_V))\geq 3r+1 -2r =r+1$$
for any $i=1,2,3$. 

From $\dim (Ker(\pi_1|_V)) \geq r+1$ we can deduce
 there would exist  in $V$ at least $r+1 $ linear independent vectors with the first $2m$ coordinates equal to $0$.

From $\dim (Ker(\pi_2|_V)) \geq r+1$ we can deduce 
 there would exist in $V$  at least $r+1 $ linear independent vectors with the last $2m$ coordinates equal to $0$.
 
But then $\dim (\pi_3(V))$ would be greater than or equal to $2(r+1)$, which is absurd.
 
\end{proof}

\begin{lem} \label{lemma3}
Let $n_1, \ldots, n_k, q_1, \ldots, q_k, m, r \in \N$.
Let $h=3m + n_1+ \ldots+ n_k$.

Let $\pi_1: \R^h \rightarrow \R^{2m}$ be the projection onto the first 
$2m$ coordinates.

Let $\pi_2: \R^h \rightarrow \R^{2m} $ be the projection onto the coordinates $m+1, \ldots, 3m$.

Let $\pi_3 :  \R^h \rightarrow \R^{2m}$ be the projection 
onto the  coordinates $1, \ldots, m, 2m+1, \ldots 3m$.

Finally, let $p_1 : \R^h \rightarrow 
\R^{n_1}$ be the projection onto the coordinates $3m+1, \ldots, 3m + n_1$, 
 let $p_2 : \R^h \rightarrow 
\R^{n_2}$ be the projection onto the coordinates $3m + n_1+1, \ldots, 3m +n_1+n_2$ and so on.

Let $V$ be a vector subspace  of $\R^h$ such that $\dim(\pi_i(V)) \leq 2r $ for $i=1,2,3$
and $\dim(p_j (V)) \leq q_j $ for $j=1, \ldots , k$; then
$$\dim(V) \leq \sum_{j=1, \ldots, k} q_j + 3r.$$
\end{lem}

\begin{proof}

Let $\tilde{\pi}_1: \R^{3m} \rightarrow \R^{2m}$ be the projection onto the first 
$2m$ coordinates.

Let $\tilde{\pi}_2: \R^{3m} \rightarrow \R^{2m} $ be the projection onto the coordinates $m+1, \ldots, 3m$.

Let $\tilde{\pi}_3 :  \R^{3m} \rightarrow \R^{2m}$ be the projection 
onto the  coordinates $1, \ldots, m, 2m+1, \ldots 3m$.

Let $\pi: \R^h \rightarrow \R^{3m}$ be the projection onto the first 
$3m$ coordinates.

Obviously $\tilde{\pi}_i \circ \pi= \pi_i$;  so, for $i=1,2,3$,
$$ \dim (\tilde{\pi}_i(\pi(V)))= \dim (\pi_i(V))
\leq 2r, $$
where the inequality holds by assumption. Hence, by Lemma \ref{lemma2}, we have that $$\dim (\pi(V)) \leq 3r .$$
Moreover 
$$ V \subseteq \pi(V)+ p_1(V)+ \ldots+ p_k(V),$$
hence 
$$ \dim (V)  \leq \dim (\pi(V)) + \dim(p_1(V))+ \ldots+ \dim(p_k(V)) \leq 3 r + q_1+ \ldots + q_k.$$

\end{proof}

\begin{lem} \label{antelemmaS}
Let $r \in \N-\{0\}$ and let $ b,c \in \R^{2r}$.
Then $$ \det \bp  \overline{J}_{2r} & c \\ {}^t b&0  \ep = c_2 b_1 - c_1 b_2 + c_4 b_3 - c_3 b_4 + \ldots $$
\end{lem}

\begin{proof}
By Schur's Complement Lemma 
$$ \det \bp  \overline{J}_{2r} & c \\ {}^t b&0  \ep = \det (\overline{J}_{2r} ) \det \left( 
- {}^t b \, \overline{J}_{2r} ^{-1} \, c 
\right)=  \det (\overline{J}_{2r} ) \det \left( 
 {}^t b \, \overline{J}_{2r}  \, c 
\right)=$$ $$= \det \left( {}^t b \bp c_2 \\ -c_1 \\ c_4\\ -c_3 \\ . \\. \\. \ep \right)
=
c_2 b_1 - c_1 b_2 + c_4 b_3 - c_3 b_4 + \ldots $$
\end{proof}

\begin{cor} \label{lemmaS}
Let $r \in \N-\{0\}$, $ b,c \in \R^{2r}-\{0\}$, $A \in M(2r \times 2r, \R)$.
Then $$ \det \bp  \overline{J}_{2r} +sA & sc \\ s {}^t b&0  \ep
$$ is a polynomial in $s$ whose lowest degree term is $$s^2(
 c_2 b_1 - c_1 b_2 + c_4 b_3 - c_3 b_4 + \ldots ).$$
\end{cor}

\begin{proof}[Proof of Theorem   
\ref{thmantisym}]

First suppose  \underline{$n \geq 2r+2$}. Consider the following affine subspace:
$$S=\left\{   
 \bp  
 \overline{J}_{2r} +R(A_{i,j}) & C \\ 
 -{}^t C & 0  \ep 
 | \; \begin{array}{l}A_{i,j} \in M(2 \times 2, \R)   \; \mbox{\rm with the second row equal to zero   } 
 \\  \;\;\;\;\; \mbox{\rm for } \; i, j \in \{1,\ldots , r\},\;   \mbox{\rm  with } i <j , \\
  C \in M(2r \times (n-2r), \R) \;\mbox{\rm   s.t.}\;C_{2i,j}=0 \\ \;\;\;\;\;  \mbox{\rm  for } i=1, \ldots , r, \; j=1, \ldots , n-2r
 \end{array}
 \right\}.$$
Observe that every element $A$ of $S $ has rank $2r$, in fact  the submatrix $$A^{(1,\ldots, 2r)}_{(1, \ldots, 2r)} =\overline{J}_{2r} +R(A_{i,j}) $$ is invertible and 
 every $(2r+1) \times (2r+1)$-submatrix  $A^{(1,	\ldots, 2r,i)}_{(1, \ldots, 2r,j )}$ for $i,j \in \{ 2r+1, \ldots ,n \}$ has determinant equal to $0$. So $S \in {\cal A}_{antisym}^{\R}(n; 2r)$.
 
  Obviously the dimension of $S$ is $r(n-2r) + r(r-1) $, i.e. 
$r (n-r-1)  $. So we get
$$a_{antisym}^{\mathbb{R}}( n; 2r) \geq (n-r-1) r .$$

If \underline{$n=2r+1$}, consider 
 the following affine subspace:
$$S=\left\{   
 \bp  
 \overline{J}_{2r} +R(A_{i,j}) & C \\ 
 -{}^t C & 0  \ep 
 | \; \begin{array}{l}A_{i,j} \in M(2 \times 2, \R)   \; \mbox{\rm with the second row equal to zero   } 
 \\  \;\;\;\;\; \mbox{\rm for } \; i, j \in \{1,\ldots , r\},\;   \mbox{\rm  with }  i <j, \\
  C \in M(2r \times 1, \R) 
 \end{array}
 \right\}.$$
Observe that every element $A$ of $S $ has rank $2r$, in fact  the submatrix $$A^{(1,\ldots, 2r)}_{(1, \ldots, 2r)} =\overline{J}_{2r} +R(A_{i,j}) $$ is invertible and obviously $A$ has determinant equal to $0$.

 The dimension of $S$ is $ r(r-1) +2r$, i.e. 
$r (r+1)  $. So, in this case,  we get
$$a_{antisym}^{\mathbb{R}}( n; 2r) \geq r(r+1) .$$

If \underline{$n=2r$}, consider 
 the following affine subspace:
$$S=\left\{   
 \bp  
 \overline{J}_{2r} +R(A_{i,j}) \\ 
  \ep 
 | \; \begin{array}{l}A_{i,j} \in M(2 \times 2, \R)   \; \mbox{\rm with the second row equal to zero   } 
 \\  \;\;\;\;\; \mbox{\rm for } \; i, j \in \{1,\ldots , r\},\;   \mbox{\rm  with }  i <j, \\
 
 \end{array}
 \right\}.$$
 Every element $A$ of $S $ has rank $2r$ and the dimension of $S$ is $ r(r-1) $. So, in this case,  we get
$$a_{antisym}^{\mathbb{R}}( n; 2r) \geq r(r-1) .$$

\medskip

Now let us prove
 the other inequalities.

Suppose $n \geq 2r+2$.

 Let $R \in  {\cal A}_{antisym}^{\mathbb{R}} (n;2r) $. 
 We want to prove that $\dim (R) \leq  r(n-r-1)
 $. We can write $R$ as
  $M+ L$ where $M \in A(n, \mathbb{R})$ and $L$ 
  is a linear subspace of $A(n, \mathbb{R})$. 
Let  $$\tilde{J} = \bp 
 \overline{J}_{2r} & 0  \\ 0 & 0 
   \ep.$$
  Let $Q$ be an invertible matrix such that ${}^t Q M Q$ is  a   matrix  of the kind $$diag(d_1, d_1, d_2, d_2, \ldots, d_r, d_r, 1, \ldots, 1) \tilde{J}$$ with $d_i \in \{-1,1\}$. 
 Let $$V= {}^t Q L Q$$ and $$S=  {}^t Q R Q = diag(d_1, d_1, d_2, d_2, \ldots, d_r, d_r, 1, \ldots, 1)\tilde{J}+V.$$ Obviously $S  \in  {\cal A}_{antisym}^{\mathbb{R}} (n;2 r) $; moreover, $ \dim(S)= \dim (R)$, so, to prove that $\dim (R) \leq  r(n-r-1) $,
  it is sufficient to prove that $$\dim (S) \leq 
 r(n-r-1).$$

Suppose for simplicity $d_1,\ldots , d_r =1$, i.e.  $S= \tilde{J}+V $. In the other cases  we can argue analogously.

\medskip

Let $$P= \{W  \in A(n,\R) | \; W_{i,j}=0 \; \mbox{\rm for} \;i,j \in 
 \{n-2r+1, \ldots, n\}\},$$
 $$K=\left\{   
 \bp  
 0_{2r \times 2r}  & C \\ 
 -{}^t C & 0_{(n-2r) \times (n-2r)}  \ep
 | \;  
  C \in M(2r \times (n-2r), \R)
 \right\}$$ 
 and, 
for any $j_1, j_2 \in \{1, \ldots,  n-2r\}$, 
 $$K_{j_1, j_2}= \left\{   
 \bp  
 0_{2r \times 2r}  & C \\ 
 -{}^t C & 0_{(n-2r) \times (n-2r)}  \ep
 | \;  
  C \in M(2r \times (n-2r), \R) \;\mbox{\rm   s.t.}\;C^{(j)}=0  \; \mbox{\rm  for } j\neq j_1, j_2
 \right\}$$
 
First let us prove that 
\begin{equation} \label{w=0}
 V\subseteq P
. \end{equation}
 Suppose that  there exists  $v \in V$ 
 such that $v_{i,j} \neq 0$ 
 for some $i,j \in \{n-2r+1, \ldots, n\}$. But in this case  
$ \det \left( (\tilde{J} + s v )^{(1,\ldots, 2r, j)}_{(1,\ldots, 2r, i)} \right)$ would be a polynomial in $s$ with coefficient of the term of degree $1$ equal to $v_{i,j} $ (which is nonzero),  so a nonconstant polynomial; thus there would exist $s$ such that  $ \det \left( (\tilde{J} + sv)^{(1,\ldots, 2r, j)}_{(1,\ldots, 2r, i)} \right) \neq 0$.
 So, for such an $s$, we would have  ${\operatorname{rk}} (\tilde{J} + s v) >2r$, 
  so $S $ 
 would not be of constant rank $2r$, which is contrary to our assumption. 
 
 Hence we  have proved (\ref{w=0}).

For any $j_1, j_2 \in \{1, \ldots,  n-2r\}$ with $j_1 \neq j_2$, 
let $Z_{j_1, j_2} $ be a $(2r)$-dimensional vector subspace of $K_{j_1, j_2}$ 
such that the quadratic form (of signature $(2r,2r)$)
\begin{equation}
\label{defpos} \sum_{i=1,\ldots, r} c_{2i-1, j_2} c_{2i ,j_1} -
c_{2i-1,j_1} c_{2i, j_2}
\end{equation}
 is positive definite 
on $Z_{j_1, j_2} $.

Let $$\pi: V \longrightarrow K$$
 be the map 
 $$  \bp  
 A  & C \\ 
 -{}^t C & 0_{(n-2r) \times (n-2r)}  \ep \longmapsto  \bp  
 0_{2r \times 2r}  & C \\ 
 -{}^t C & 0_{(n-2r) \times (n-2r)}  \ep
 $$

For any $j_1, j_2 \in \{1, \ldots,  n-2r\}$ with $j_1 \neq j_2$,  let 
$$\pi_{j_1, j_2}: K
\longrightarrow 
K_{j_1, j_2}$$ be the map sending every matrix 
$  \bp  
 0& C \\ 
 -{}^t C & 0_{(n-2r) \times (n-2r)}  \ep$ in $K$  to  the matrix
$$\bp  
 0_{2r \times 2r } & \left( 0\ldots0, C^{(j_1)},0 \ldots 0, C^{(j_2)} ,0 \ldots 0 \right) \\ 
- {}^t \left(  0\ldots0, C^{(j_1)},0 \ldots 0, C^{(j_2)} ,0 \ldots 0 \right) & 0_{(n-2r) \times (n-2r)}  \ep .$$

 Fix $j_1, j_2 \in \{1, \ldots,  n-2r\}$ with $j_1 \neq j_2$. 
 We want to prove that 
 \begin{equation} \label{zj1j2} \pi_{j_1, j_2}(
 \pi(V)) \cap Z_{j_1, j_2}=\{0\}
\end{equation}

Suppose on the contrary that there exists 
 $x=\bp  
 A  & C \\ 
 -{}^t C & 0_{(n-2r) \times (n-2r)}  \ep  \in V$ such that $\pi_{j_1, j_2}( \pi(x)) \in  Z_{j_1, j_2} \setminus \{0\}$.
 Then the number in (\ref{defpos})
 would be positive. But, in this case, we would have,  by Corollary \ref{lemmaS}, that 
$$ \det \left( (\tilde{J} + sx )^{(1,\ldots, 2r, j_2+2r)}_{(1,\ldots, 2r, j_1+2r)} \right)$$ is  
a polynomial in $s$ with coefficient of the term of degree $2$ equal to the number in (\ref{defpos})
 (which is positive), hence a nonconstant polynomial.
Thus there would exist $s$ such that  
$ \det \left( (\tilde{J} + s x)^{(1,\ldots, 2r, j_2+2r)}_{(1,\ldots, 2r, j_1+2r)} \right) \neq 0$ 

 So, for such a $s$,  we would have ${\operatorname{rk}} (\tilde{J} + sx) > 2r$, 
  so $S $ 
 would not be of constant rank $2r$, which is contrary to our assumption.

 Hence we have proved (\ref{zj1j2}); therefore
\begin{equation} \label{dim}
\dim 
( \pi_{j_1, j_2} (\pi(V))) \leq 2r .
\end{equation}

If $n-2r $ is even consider the projections
$\pi_{1,2}, \pi_{3,4}, \ldots, \pi_{n-2r-1, n-2r}$.
If $n-2r $ is odd consider the projections
$\pi_{1,2}, \pi_{1,3}, \pi_{2,3}, \pi_{4,5}, \ldots, \pi_{n-2r-1, n-2r}$.

By Lemma \ref{lemma3} and from (\ref{dim})
we get that $$\dim(\pi(V)) \leq r (n-2r).$$
Hence there exists a  $ r (n-2r)$-dimensional
vector subspace $Z $ in $K$
such that 
 \begin{equation} \label{Z}
 \pi(V) \cap Z =\{0\}.
 \end{equation}

 Let $$U=
 \left\{   
 \bp  
 T(l_1, \ldots, l_{r} ) +R(A_{i,j}) & 0_{2r \times (n-2r)}  \\ 
 0_{(n-2r) \times 2r}  & 0_{(n-2r) \times (n-2r)}  \ep 
 | \; \begin{array}{l}A_{i,j} \; 2 \times 2  \; \mbox{\rm  antipinco matrix} \\ 
\mbox{\rm    for} \; i, j \in \{1,\ldots , r\},\; i <j,  \\l_1, \ldots, l_{r} \in \R
 \end{array}
 \right\}.$$
Observe that $$\dim(U)= r(r-1)+r= r^2.$$
  
  We want to prove that $$ V \cap (U+Z)=\{0\}$$

 Let $u \in U $ and $z \in Z$ be such that $u+z \in V$. 
 
 The element  $z$ must be zero by formula (\ref{Z}).
 
Finally suppose $u \neq 0$.  
Let $u= 
 \bp  
 T(l_1, \ldots, l_{2r-1} ) +R(A_{i,j}) & 0_{2r \times (n-2r)}  \\ 
 0_{(n-2r) \times 2r}  & 0_{(n-2r) \times (n-2r)}  \ep $ with 
 $A_{i,j}$  $ 2 \times 2 $ antipinco matrices, $l_1,\ldots, l_r \in \R$.
 
 Observe that
  
 $ J \bp a & b \\ -b & a \ep = \bp -b & a \\ -a & -b\ep $,

 $ J \left( -{}^t \bp a & b \\ -b & a \ep \right) = \bp -b & -a \\ a & -b\ep \  $,

$J \bp 0 & l \\ -l & 0 \ep = -l I$.

Hence
\begin{equation} \label{ss}
-\overline{J}_{2r} \left( 
\overline{J}_{2r}+ s [ T(l_1, \ldots, l_{r} ) +R(A_{i,j})]
\right)= I_{2r} + s A,
\end{equation}
 where $$ A
= diag (l_1, l_1, l_2, l_2, \ldots, l_r , l_r) + X(T_{i,j}) $$ for some  antipinco matrices $T_{i,j}$ for $i, j \in \{1, \ldots, r\}$ with $i<j$.
  
Since $A$ is symmetric and nonzero (since $u \neq 0$), it has an eigenvalue $\lambda \in \R-\{0\}$; let $s = - \frac{1}{\lambda}$; we have:
$$ \det(I_{2r} + s A)= s^{2r} \left(\det \left(A+
 \frac{1}{s}I_{2r} \right)\right)= s^{2r} \left(\det \left(A- \lambda
 I_{2r} \right)\right)=0. $$
 But $ \det(I_{2r} + s A)= \det  \left( 
\overline{J}_{2r}+ s \left[T(l_1, \ldots, l_{r} ) +R(A_{i,j}) \right]
\right)$ by (\ref{ss}); hence we get $$ \det  \left( 
\overline{J}_{2r}+ s \left[T(l_1, \ldots, l_{r} ) +R(A_{i,j}) \right]
\right)=0$$
and  then ${\operatorname{rk}} (\tilde{J} + s(u+z) )=  {\operatorname{rk}} (\tilde{J} + su ) < 2r$, 
  so $S $ 
 would not be of constant rank $2r$, which is contrary to our assumption. Therefore also $u$ must be equal to $0$.

So we have proved that 
$V \cap (U + Z)=\{0\}$. 

Since $V$ and $U+Z$ are subspaces of $P$, we can deduce that $\dim (V)$ is less than or equal to 
 $$ \dim (P) -\dim (U+Z)=$$ 
  $$ =  \dim (P)  -\dim (U) - \dim(Z)= 
 $$ 
 $$= r(2r-1) + 2r (n-2r)  - r^2 - r(n-2r)=  r(n-r-1)$$
 
 In the cases $n=2r$ and $n=2r+1$, we can argue analogously (in this case we have only to show that $V \cap U=\{0\}$ and then to use that $\dim(V ) \leq \dim (P) - \dim (U)$).
 
\end{proof}

{\bf Acknowledgments.}
This work was supported by the National Group for Algebraic and Geometric Structures, and their  Applications (GNSAGA-INdAM).


\end{document}